\documentclass[10pt]{article}
\usepackage{amssymb}
\usepackage{amsfonts}
\usepackage{color}

\newcommand{\nc}{\newcommand}
\nc{\thref}[1]{Theorem~\ref{theo:#1}}
\nc{\selabel}[1]{\label{sect:#1}}
\nc{\seref}[1]{Section~\ref{sect:#1}}
\nc{\lelabel}[1]{\label{lemm:#1}}
\nc{\leref}[1]{Lemma~\ref{lemm:#1}}
\nc{\prlabel}[1]{\label{prop:#1}}
\nc{\prref}[1]{Proposition~\ref{prop:#1}}
\nc{\colabel}[1]{\label{coro:#1}}
\nc{\coref}[1]{Corollary~\ref{coro:#1}}
\nc{\exlabel}[1]{\label{exam:#1}}
\nc{\exref}[1]{Example~\ref{exam:#1}}
\nc{\delabel}[1]{\label{defi:#1}}
\nc{\deref}[1]{Definition~\ref{defi:#1}}
\nc{\eqlabel}[1]{\label{equa:#1}}
\nc{\relabel}[1]{\label{rema:#1}}
\nc{\reref}[1]{Lemma~\ref{rema:#1}}
\providecommand{\operatorname}[1]{\mathrm{#1}\,}
\nc{\Hom}{\operatorname{Hom}}
\nc{\Mor}{\operatorname{Mor}}
\nc{\Aut}{\operatorname{Aut}}
\nc{\Ann}{\operatorname{Ann}}
\nc{\Ker}{\operatorname{Ker}}
\nc{\Trace}{\operatorname{Trace}}
\nc{\Char}{\operatorname{Char}}
\nc{\Mod}{\operatorname{Mod}}
\nc{\End}{\operatorname{End}}
\nc{\Spec}{\operatorname{Spec}}
\nc{\Span}{\operatorname{Span}}
\nc{\sgn}{\operatorname{sgn}}
\nc{\Id}{\operatorname{Id}}
\nc{\Com}{\operatorname{Com}}
\nc{\rank}{\operatorname{rank}}

\let\:=\colon

\newtheorem{de}{Definition}[section]
\newtheorem{lm}[de]{Lemma}
\newtheorem{pr}[de]{Proposition}
\newtheorem{co}[de]{Corollary}
\newtheorem{re}[de]{Remark}
\newtheorem{res}[de]{Remarks}
\newtheorem{te}[de]{Theorem}
\newtheorem{ex}[de]{Example}
\newtheorem{exs}[de]{Examples}

\def\bex{\begin{ex}}
\def\eex{\end{ex}}
\def\bexs{\begin{exs}}
\def\eexs{\end{exs}}
\def\bl{\begin{lm}}
\def\el{\end{lm}}
\def\bc{\begin{co}}
\def\ec{\end{co}}
\def\bt{\begin{te}}
\def\et{\end{te}}
\def\bpr{\begin{pr}}
\def\epr{\end{pr}}
\def\br{\begin{re}}
\def\er{\end{re}}
\def\brs{\begin{res}}
\def\ers{\end{res}}
\def\bd{\begin{de}}
\def\ed{\end{de}}
\def\be{\begin{equation}}
\def\ee{\end{equation}}
\def\bea{\begin{eqnarray*}}
\def\eea{\end{eqnarray*}}
\def\bp{\begin{proof}}
\def\ep{\end{proof}}

\def\qed{\hfill\Box}

\let\:=\colon

\title{\bf{A new and elementary proof of Newton's "favorite" quadrature formulae}}
\author{\bf{Cezar Lupu, Tudorel Lupu}}
\date{}

\begin{document}
\maketitle

In this note we shall give a new proof to a quadrature formulae due to Newton(see [1.]) which states

\bt\label{t1}
If $f:[a, b]\to\mathbb{R}$ is four times differentiable on $[a, b]$ with $f^{(4)}$ continuous, then there is some $\xi\in [a,b]$ such that
$$\int_{a}^{b}f(x)dx=\frac{b-a}{8}\left[f(a)+3f\left(\frac{2a+b}{3}\right)+3f\left(\frac{2b+a}{3}\right)+f(b)\right]-\frac{(b-a)^{5}}{6480}f^{(4)}(\xi).$$

\et

In order to prove this theorem we need to show first three lemmas starting with\\

\noindent {\bf{Lemma 1.}}
If $f:[a,b]\to\mathbb{R}$ is four times differentiable with $f^{(4)}$ continuous, then there is some $\xi_{1}\in [a,b]$ such that
$$\int_{a}^{b}f(x)dx=(b-a)f\left(\frac{a+b}{2}\right)+\frac{(b-a)^{3}}{24}f^{(2)}\left(\frac{a+b}{2}\right)+\frac{(b-a)^{5}}{1920}f^{(4)}(\xi_{1}).$$

\noindent {\bf{Proof.}}
We begin with the following integral
$$\int_{a}^{b} (x-a)^{4}\left[f^{(4)}\left(\frac{a+x}{2}\right)+f^{(4)}\left(\frac{2b+a-x}{2}\right)\right]dx.$$
Integrating by part we get
$$\int_{a}^{b} (x-a)^{4}\left[f^{(4)}\left(\frac{a+x}{2}\right)+f^{(4)}\left(\frac{2b+a-x}{2}\right)\right]dx=2(x-a)^{4}\left[f^{(3)}\left(\frac{a+x}{2}\right)-f^{(3)}\left(\frac{2b+a-x}{2}\right)\right]dx\arrowvert_a^b$$
$$-2\int_{a}^{b}4(x-a)^{3}\left[f^{(3)}\left(\frac{a+x}{2}\right)-f^{(3)}\left(\frac{2b+a-x}{2}\right)\right]=-16(x-a)^{3}\left[f^{(2)}\left(\frac{a+x}{2}\right)+f^{(2)}\left(\frac{2b+a-x}{2}\right)\right]\arrowvert_a^b$$
$$+48\int_{a}^{b}(x-a)^{2}\left[f^{(2)}\left(\frac{a+x}{2}\right)+f^{(2)}\left(\frac{2b+a-x}{2}\right)\right]dx=-32(b-a)^{3}f^{(2)}\left(\frac{a+b}{2}\right)+$$
$$+96(x-a)^{2}\left[f'\left(\frac{a+x}{2}\right)-f'\left(\frac{2b+a-x}{2}\right)\right]\arrowvert_a^b-192\int_{a}^{b}(x-a)\left[f'\left(\frac{a+x}{2}\right)-f'\left(\frac{2b+a-x}{2}\right)\right]dx$$
$$=-32(b-a)^{3}f^{(2)}\left(\frac{a+b}{2}\right)-384(x-a)\left[f\left(\frac{a+x}{2}\right)+f\left(\frac{2b+a-x}{2}\right)\right]\arrowvert_a^b$$
$$+384\int_{a}^{b}\left[f\left(\frac{a+x}{2}\right)+f\left(\frac{2b+a-x}{2}\right)\right]dx=768\int_{a}^{b} f(x)dx-768(b-a)f\left(\frac{a+b}{2}\right)-32(b-a)^{3}f^{(2)}\left(\frac{a+b}{2}\right).$$
On the other hand, by the mean value theorem, there is $\eta\in [a, b]$ such that
$$\int_{a}^{b}(x-a)^{4}\left[f^{(4)}\left(\frac{a+x}{2}\right)+f\left(\frac{2b+a-x}{2}\right)\right]dx=\left[f^{(4)}\left(\frac{a+\eta}{2}\right)+f^{(4)}\left(\frac{2b+a-\eta}{2}\right)\right]\int_{a}^{b}(x-a)^{4}dx.$$
But, $f^{(4)}$ is continuous, so $f^{(4)}$ has intermediate value property. It follows that there is $\xi_{1}\in [a,b]$ such that
$$f^{(4)}\left(\frac{a+\eta}{2}\right)+f^{(4)}\left(\frac{2b+a-\eta}{2}\right)=2f^{(4)}(\xi_{1}),$$
so $$\int_{a}^{b}(x-a)^{4}\left[f^{(4)}\left(\frac{a+x}{2}\right)+f\left(\frac{2b+a-x}{2}\right)\right]dx=\frac{2(b-a)^{5}}{5}f^{(4)}(\xi_{1}).$$
We get
$$768\int_{a}^{b}f(x)dx=768(b-a)f\left(\frac{a+b}{2}\right)+32(b-a)^{3}f^{(2)}\left(\frac{a+b}{2}\right)+\frac{2(b-a)^{5}}{5}f^{(4)}(\xi_{1})$$
which is equivalent to
$$\int_{a}^{b}f(x)dx=(b-a)f\left(\frac{a+b}{2}\right)+\frac{(b-a)^{3}}{24}f^{(2)}\left(\frac{a+b}{2}\right)+\frac{(b-a)^{5}}{1920}f^{(4)}(\xi_{1}),$$
where $\xi_{1}\in [a,b]$.     $\qed$

\noindent {\bf{Lemma 2.}}

If $f:[a,b]\to\mathbb{R}$ is four times differentiable with $f^{(4)}$ continuous, then there is some $\xi_{2}\in [a,b]$ such that
$$\int_{a}^{b}f(x)dx=\frac{(b-a)}{8}\left[f(a)+6f\left(\frac{a+b}{2}\right)+f(b)\right]+\frac{(b-a)^{3}}{96}f^{(2)}\left(\frac{a+b}{2}\right)-\frac{(b-a)^{5}}{7680}f^{(4)}(\xi_{2}).$$

\noindent {\bf{Proof.}}
Like in the proof of lemma 1, we will begin with the following integral
$$\int_{a}^{b}(x-a)^{3}(x-b)\left[f^{(4)}\left(\frac{a+x}{2}\right)+f^{(4)}\left(\frac{2b+a-x}{2}\right)\right]dx.$$
Integrating by part we get
$$\int_{a}^{b}(x-a)^{3}(x-b)\left[f^{(4)}\left(\frac{a+x}{2}\right)+f^{(4)}\left(\frac{2b+a-x}{2}\right)\right]dx=2(x-a)^{3}(x-b)\left[f^{(3)}\left(\frac{a+x}{2}\right)-f^{(3)}\left(\frac{2b+a-x}{2}\right)\right]\arrowvert_a^b-$$
$$-2\int_a^b\left[3(x-a)^{2}(x-b)+(x-a)^{3}\right]\left[f^{(3)}\left(\frac{a+x}{2}\right)-f^{(3)}\left(\frac{2b+a-x}{2}\right)\right]dx=$$
$$=-4\left[3(x-a)^{2}(x-b)+(x-a)^{3}\right]\left[f^{(2)}\left(\frac{a+x}{2}\right)+f^{(2)}\left(\frac{2b+a-x}{2}\right)\right]\arrowvert_a^b+$$
$$+24\int_{a}^{b}\left[(x-a)^{2}+(x-a)(x-b)\right]\left[f^{(2)}\left(\frac{a+x}{2}\right)+f^{(2)}\left(\frac{2b+a-x}{2}\right)\right]dx=$$
$$=-8(b-a)^{3}f^{(2)}\left(\frac{a+b}{2}\right)+48\left[(x-a)^{2}+(x-a)(x-b)\right]\left[f'\left(\frac{a+x}{2}\right)-f'\left(\frac{2b+a-x}{2}\right)\right]\arrowvert_a^b-$$
$$-48\int_{a}^{b}\left[3(x-a)+x-b\right]\left[f'\left(\frac{a+x}{2}\right)-f'\left(\frac{2b+a-x}{2}\right)\right]dx=$$
$$=-8(b-a)^{3}f^{(2)}\left(\frac{a+b}{2}\right)-96\left[3(x-a)+(x-b)\right]\left[f\left(\frac{a+x}{2}\right)+f\left(\frac{2b+a-x}{2}\right)\right]\arrowvert_a^b$$
$$+384\int_{a}^{b}\left[f\left(\frac{a+x}{2}\right)+f\left(\frac{2b+a-x}{2}\right)\right]dx=768\int_{a}^{b}f(x)dx-96\cdot 6(b-a)f\left(\frac{a+b}{2}\right)$$
$$-96(b-a)\left[f(a)+f(b)\right]-8(b-a)^{3}f^{(2)}\left(\frac{a+b}{2}\right)=768\int_{a}^{b}f(x)dx-96(b-a)\left[f(a)+6f\left(\frac{a+b}{2}\right)+f(b)\right]-$$
$\displaystyle -8(b-a)^{3}f^{(2)}\left(\frac{a+b}{2}\right)$.
On the other hand, by the mean value theorem there is some $\eta'\in [a,b]$ such that
$$\int_{a}^{b}(x-a)^{3}(x-b)\left[f^{(4)}\left(\frac{a+x}{2}\right)+f^{(4)}\left(\frac{2b+a-x}{2}\right)\right]dx=\left[f^{(4)}\left(\frac{a+\eta'}{2}\right)+f^{(4)}\left(\frac{2b+a-\eta'}{2}\right)\right]$$
$\displaystyle\cdot\int_{a}^{b}(x-a)^{3}(x-b)dx$.
By the intermediate value property, there is some $\xi_{2}\in [a,b]$ such that
$$f^{(4)}\left(\frac{a+\eta'}{2}\right)+f^{(4)}\left(\frac{2b+a-\eta'}{2}\right)=2f^{(4)}(\xi_{2}).$$
A simple calculation shows that $\displaystyle\int_{a}^{b}(x-a)^{3}(x-b)dx=-\frac{(b-a)^{5}}{20}$.
Now, combining these two equalities obtained, the lemma follows immediately.      $\qed$
\\

\noindent {\bf{Lemma 3}(Simpson-Cavalieri).}

If $f:[a,b]\to\mathbb{R}$ is four times differentiable with $f^{(4)}$ continuous, then there is some $\xi_{3}\in [a,b]$ such that
$$\int_{a}^{b}f(x)dx=\frac{(b-a)}{6}\left[f(a)+4f\left(\frac{a+b}{2}\right)+f(b)\right]-\frac{(b-a)^{5}}{2880}f^{(4)}(\xi_{3}).$$

\noindent {\bf{Proof.}}
This formulae can be proved in many ways, but we prefer to deduce it by eleminating $\displaystyle f^{(2)}\left(\frac{a+b}{2}\right)$ from lemma 1 and lemma 2.    $\qed$

Now, we shall prove the main result of our paper. First, we apply lemma 1 on the interval $\displaystyle\left[\frac{2a+b}{3},\frac{2b+a}{3}\right]$.
Then, there is some $\displaystyle \xi'_{1}\in \left[\frac{2a+b}{3},\frac{2b+a}{3}\right] $ such that
$$\int_{\frac{2a+b}{3}}^{\frac{2b+a}{3}}f(x)dx=\frac{b-a}{3}f\left(\frac{a+b}{2}\right)+\frac{(b-a)^{3}}{24\cdot 27}f^{(2)}\left(\frac{a+b}{2}\right)+\frac{(b-a)^{5}}{1920\cdot     3^5}f^{(4)}(\xi'_{1}).(*)$$
By eliminating $\displaystyle f^{(2)}\left(\frac{a+b}{2}\right)$ from lemma 2 and $(*)$ we get
$$4\int_{a}^{b}f(x)dx=27\int_{\frac{2a+b}{3}}^{\frac{2b+a}{3}}f(x)dx+\frac{b-a}{2}\left[f(a)-12f\left(\frac{a+b}{2}\right)+f(b)\right]-\frac{(b-a)^{5}}{9\cdot 192}f^{(4)}(\xi_{4}),(**)$$
with $\xi_{4}\in [a,b]$.
Now, applying Lemma 3(Simpson-Cavalieri), we get 
$$\int_{\frac{2a+b}{3}}^{\frac{2b+a}{3}}f(x)dx=\frac{b-a}{18}\left[f\left(\frac{2a+b}{3}\right)+4f\left(\frac{a+b}{2}\right)+f\left(\frac{2b+a}{3}\right)\right]-\frac{(b-a)^{5}}{3^5\cdot 2880}f^{(4)}(\xi_{5})(***)$$
with $\displaystyle\xi_{5}\in\left[\frac{2a+b}{3},\frac{2b+a}{3}\right]$. 
Now, replacing $(***)$ in $(**)$ we have
$$4\int_{a}^{b}f(x)dx=\frac{27(b-a)}{18}\left[f\left(\frac{2a+b}{3}\right)+4f\left(\frac{a+b}{2}\right)+f\left(\frac{2b+a}{3}\right)\right]-\frac{(b-a)^{5}}{9\cdot 2880}f^{(4)}(\xi_{5})$$
$$+\frac{(b-a)}{2}\left[f(a)-12f\left(\frac{a+b}{2}\right)+f(b)\right]-\frac{(b-a)^{5}}{9\cdot 192}f^{(4)}(\xi_{4}),$$
with $\xi_{4}, \xi_{5}\in [a,b]$.
The above formulae is equivalent to
$$4\int_{a}^{b}f(x)dx=\frac{b-a}{2}\left[f(a)+3f\left(\frac{2a+b}{3}\right)+3f\left(\frac{2b+a}{3}\right)+f(b)\right]-\left[\frac{(b-a)^{5}}{9\cdot 2880}f^{(4)}(\xi_{5})+\frac{(b-a)^{5}}{9\cdot 192}f^{(4)}(\xi_{4})\right].$$
Now, using intermediate value property, we get a $\xi\in [a,b]$ such that
$$\frac{(b-a)^{5}}{9\cdot 2880}f^{(4)}(\xi_{5})+\frac{(b-a)^{5}}{9\cdot 192}f^{(4)}(\xi_{4})=\frac{(b-a)^{5}}{1620}f^{(4)}(\xi).$$
Now, Newton's result follows immediately.

\noindent {\bf{References.}}

[1.] D.V. Ionescu-Restul \^ in formula de cuadratur\u a "preferat\u a" a lui Newton, Gaz.Mat., vol 7, 1965, pag. 298-303.

[2.] J.V. Steffensen-Interpolation, Baltimore, 1927

[3.] B.P. Demidovich, I.A Maron-Computational Mathematics, Moscow: Mir Publishers, 1976. 2nd edition. Translated from Russian.\\

{\bf{Cezar Lupu}, student Faculty of Mathematics, University of Bucharest, Bucharest, Romania.}
{\bf e-mail: lupucezar@yahoo.com}

{\bf{Tudorel Lupu}, teacher, Decebal High school, Constanza, Romania.}

\end{document}